\documentclass[12pt]{article}
\usepackage{amssymb,a4,latexsym}

\newtheorem{satz}{Satz}[section]
\newtheorem{Thm}[satz]{Theorem}
\newtheorem{Rem}[satz]{Remark}
\newtheorem{Lemma}[satz]{Lemma}
\newtheorem{Defi}[satz]{Definition}
\newenvironment{Bew}{{\bf \noindent Proof:}}
       {\hfill{{$\Box$}}}
\newcommand{\N}{{\mathbb N}}
\newcommand{\Z}{{\mathbb Z}}
\newcommand{\R}{{\mathbb R}}
\newcommand{\nkN}{n_k^N}
\newcommand{\dnachd}[1]{\frac{\partial}{\partial{#1}}}
\newcommand{\Skalprod}[2]{\langle #1 , #2\rangle}

\begin{document}

\title{A Markov jump process approximation of the
stochastic Burgers equation\footnote{This work was partly
supported by the NSF Grants DMS-0209326 and DMS-0139073. }}

\author{ {\it Christoph Gugg}\\
         IADM, Fachbereich Mathematik\\
         Universit\"at Stuttgart\\
         70550 Stuttgart, Germany\\
         {\tt christoph.gugg@gmx.de}\\
   \and {\it Jinqiao Duan}\\
        Department of Applied Mathematics \\
        Illinois Institute of Technology \\
        Chicago, IL 60616, USA \\
        {\tt duan@iit.edu}}
\maketitle

\abstract{ {\bf Stochastics and Dynamics}, 4(2004),245--264.

We consider the stochastic Burgers equation
\begin{equation}\label{absburg}
\dnachd{t} \psi(t,r) = \Delta \psi(t,r) +
\nabla \psi^2(t,r)+\sqrt{\gamma\psi(t,r)} \eta(t,r)
\end{equation}
with periodic boundary conditions, where $t \ge 0,$ $r \in [0,1],$
and $\eta$ is some space-time white noise. A certain
Markov jump process is constructed to approximate a solution of
this equation.}

\section{Introduction}

Scientific and engineering systems are often subject to
uncertainty or random influence. Randomness can have delicate
impact on the overall evolution of such systems.
Taking stochastic effects
into account is of central importance for the development of
mathematical models of complex phenomena in engineering and
science.  Macroscopic models in the form of  partial differential
equations for these systems contain such randomness as stochastic
forcing, uncertain parameters, random sources or inputs, and
random   boundary conditions. Stochastic partial differential
equations (SPDEs) are appropriate models for randomly influenced
systems.

Most of stochastic partial differential equations models are
nonlinear in nature. Especially the Burgers equation with
stochastic noise has attracted considerable attention, for example,
as a simplified model of fluid turbulence. Due to the nonlinearity,
numerical simulations are often necessary in order to understand
the dynamical behavior of the stochastic models. In this paper, we
propose a Markov chain approximation method for a stochastic
Burgers equation and prove its convergence.

More explicitly, our reaction-diffusion model is constructed by
dividing the unit interval into $N$ cells of length $1/N.$
We place an initial distribution of approximately $N l$ particles
into the cells. The particles in each cell independently jump to
neighboring cells according to Poisson processes with rates
$N^2 n_k$ where $n_k$ is the number of particles in cell $k,$ and
are born or die with rates $\gamma N l n_k/2.$ Moreover, to
obtain the desired nonlinearity we allow particles to jump to the
cell next to them on the left-hand side at a rate approximately
given by $N n_k^2/l.$ Our approximating process is given by a
step-function-valued process $X^N$ defined by the rescaled
``densities'' $n_k/l.$ We mainly assume that $l \ge cN.$ Then we
show that for $N \to \infty$ there exists a limit satisfying
(\ref{absburg}).

Our method bases on a work by D. Blount \cite{Blount96}. He
obtained a process solving the SPDE
\begin{equation}\label{introblount}
\dnachd{t} \psi(t,r) = \Delta \psi(t,r)
- d \psi^2(t,r)+\alpha \psi(t,r)
+ \sqrt{\psi(t,r)} \eta(t,r).
\end{equation}
(where $d\ge 0$ and $\eta$ is some space-time white noise) as a
high-density limit of a Markov jump process
consisting of birth- and death-processes and diffusion processes
similar to the jump process described above.
We verify some important martingale relationships between the
approximating Markov jump process and its generator by the method
of \cite{Ethier}. This allows writing the process approximately
as
$$
X^N(t)-X^N(0) = \int_0^t \Delta_N X^N(s) + \nabla^+_N (X^N(s))^2
d s + Z^N(t)
$$
where $Z^N$ is a mean-0-martingale and $\Delta_N$ and $\nabla^+_N$
are a discretized Laplace operator and a discretized first
derivative, respectively. It turns out that $Z^N$ consists of a
part originating from diffusion and a part coming from the birth
process, where the diffusion part vanishes in the limit. By the
method of \cite{Blount96} we can show tightness of the reaction
part of $Z^N$ in spaces $D(0,T,H^\alpha(0,1))$ with Skorohod
metric, where $H^\alpha(0,1)$ are certain Sobolev spaces. To show
tightness of the remaining part of the approximating process $X^N$
in $L^2(0,T,L^2(0,1))$ we adapt a method of \cite{Cap} and
\cite{Vishik} and especially show a discretized version of the
compactness result \cite{Gatarek}. The representation of the limit,
which is now inferred from the theorems of Prokhorov and Skorohod,
as solutions of (\ref{absburg}) follows by an application of the
theory of super Brownian motion, see \cite{Daw}, \cite{Walsh}, and
\cite{Konno}.

Markov jump process approximations of reaction-diffusion equations
have been studied for a long time.
A deterministic reaction-diffusion equation
with polynomial nonlinearities is treated in \cite{Arnold}.
The approximation of a linear reaction-diffusion equation by
space-time jump Markov processes is investigated by D. Blount and
P. Kotelenez e.g. in \cite{Kotelenez86a}, \cite{Kotelenez86b}, and
\cite{Blount91}, for various assumptions on the initial density of
particles and the number of cells,
and in different function spaces, and central limit
theorems are proved. These results are generalized to
reaction-diffusion equations with polynomial nonlinearities in
\cite{Kotelenez88}, \cite{Blount92}, and \cite{Blount93} by these
authors. In \cite{Blount94} laws of large numbers in a high density
and in a low density limit and a central limit theorem is given for
Equation (\ref{introblount}) without noise.
Only recently, M. Kouritzin and H. Long \cite{Kouritzin02} generalized
the ansatz to a much broader class of nonlinearities and applied the
idea to a reaction-diffusion equation that is driven by a Poisson
point process and describes water pollution. However, their
nonlinearities do not involve spatial derivatives. Our work seems to
be the first step in this direction.

When our work was almost finished,
we learned about a preprint by G. Bonnet
and R. Adler, \cite{Bonnet}, where Equation (\ref{absburg}) is
studied on the entire real line. Their approach is based on a
multidimensional stochastic differential equation driven by
(multiplicative) white-in-time noise. By means of Green function
representation and a tightness argument convergence of a
subsequence of solutions of the approximating SDE towards a
solution of (\ref{absburg}) is shown.

Moreover, the classical Burgers equation has been investigated in
the probability literature in a number of ways, e.g. as limit of
an asymmetric simple exclusion process or as limit of certain
particle systems driven by Brownian motions. We cannot give a
complete survey on the vast literature in this field. See e.g.
\cite{Bertini}, \cite{Calderoni}, \cite{DeMasi89}, \cite{DeMasi},
\cite{Kipnis}, \cite{Liggett}, and \cite{Oel85}, just to name a
few. An approximation of the 2-D-Navier-Stokes equation is found
in \cite{Meleard}.

Our work is organized as follows. In Section
\ref{eins}, we construct the Markov chain approximations to
(\ref{absburg}) in the manner of the above mentioned works. Section
\ref{zwei} contains the proofs of these results and in Section
\ref{drei} we establish some auxiliary results.

\section{Problem and Result}\label{eins}

In this section we introduce our models and present the main result.
\bigskip

{\bf The stochastic model:} is the stochastic Burgers
equation

\begin{eqnarray}\label{stochburgeq}
\nonumber&&\dnachd{t} \psi(t,r) = \Delta \psi(t,r) +
\nabla \psi^2(t,r)+ \sqrt{\gamma\psi(t,r)} \eta(t,r),\\
&&\psi:[0,T]\times [0,1] \to \R,
\end{eqnarray}

with initial condition $\psi(0,r)=\psi_0(r)$ and periodic
boundary conditions. $\eta$ is some space-time-white noise, and
$\Delta$ and $\nabla$ denote
$\frac{\partial^2}{\partial r^2}$ and $\dnachd{r},$ respectively.

\bigskip

{\bf The approximation model:} is a Markov jump process defined as
follows. Divide $[0,1]$ into $N$ cells of width $1/N.$ $[0,1]$ is
from now on is identified with a circle of circumference $1,$ to
obtain periodic boundary conditions. We place an initial
distribution of approximately $N l$ particles into the cells,
corresponding to the initial conditions given in the sequel, so $l$ can
be seen as initial average number of particles in a cell. For
$1 \le k\le N$ and $t \ge 0$ let $\nkN(t)$ be the number of particles
in cell $k$ at time $t.$ We suppress the $l$-dependence of $\nkN$
in our notation. Let $n^N(t) = (n_1^N(t),\ldots,n_N^N(t))$ in
$\N_0^N$. ($\N_0=\{0,1,2,3,\ldots\}=\N\cup\{0\}.$)
Define the jump rates for $n^N(t)$ by

\begin{eqnarray}\label{jumprates}
\nonumber
(n_{k-1},n_k) &\to& (n_{k-1}+1,n_k-1)
\\
\nonumber
&&\mbox{ at rate }
N^2 n_k + \frac{N}{3 l} (n_k^2+n_k n_{k-1}+n_{k-1}^2),
\\
(n_{k},n_{k+1}) &\to& (n_k-1,n_{k+1}+1)
\quad\mbox{ at rate } N^2 n_k.
\\
\nonumber
n_{k} &\to& n_{k}+1 \quad\quad\quad\quad\quad\quad
\mbox{ at rate } \gamma N l n_k/2.
\\
\nonumber
n_{k} &\to& n_{k}-1 \quad\quad\quad\quad\quad\quad
\mbox{ at rate } \gamma N l n_k/2.
\end{eqnarray}

(Observe the periodic boundary conditions for $n^N,$ i.e.
$n_{k+zN} = n_k,$ $z,k \in \Z.$)
For an introduction in Markov jump processes see, e.g.,
\cite{Ethier}. The state space of the
process is $E:=\N_0^N.$

The generator of the process $n^N(t)$ is given by

\begin{equation}\label{generator}
A(i,j) = \left\{\begin{array}{l@{\quad}l} -\lambda(i) & i=j \\
\lambda(i)Q(i,j) & i \neq j, \end{array}\right.
\end{equation}

where $i, j$ are elements of the state space $E.$ Let
$i=(n_1,\ldots,n_N),$
then $\lambda(i)$ is the sum over the rates in (\ref{jumprates}),
$\lambda(i)=\sum_{k=1}^N
N^2 n_k + \frac{N}{3 l} (n_k^2+n_k n_{k-1}+n_{k-1}^2) + N^2 n_k
+ \gamma N l n_k.$
The time the process remains in state $i$ until the next jump is
exponentially distributed with parameter $\lambda(i).$ $Q(i,j)$ is
the transition function of the underlying Markov chain
corresponding to the states of the process. If state $j$ can be
reached from state $i$, then
$Q(i,j)=\mbox{rate}(i,j)/\lambda(i)$, otherwise $Q(i,j)=0.$ If,
for instance, $j=i$ up to a jump of one particle
from a cell to a neighboring cell, that means
$j-i=(0,..,1,-1,0..)$ for instance,
where the $-1$ is at position $k$ then
$\mbox{rate}(i,j)= N^2 n_k
+ \frac{N}{3 l} (n_k^2+n_k n_{k-1}+n_{k-1}^2).$ By

\begin{equation}\label{generator2}
A f (i) = \sum_{j \in E} (f(j)-f(i)) A(i,j)
\end{equation}

$A$ operates on the real valued functions $f:E\to\R,$ see
\cite{Ethier}.  From \cite{Ethier}, Prop. 4.1.7, e.g., we obtain
that $f(n^N(t)) - \int_0^t A f (n^N(s)) d s$ is a martingale
w.r.t.
the filtration ${\cal F}^N_t \subset {\cal F}$ on the underlying
probability space $(\Omega,{\cal F}, P)$ which is the completion
of the $\sigma$-field induced by the
process $n^N(t).$ Let $f=f_k,$ $f_k(n_1,\ldots,n_N)=n_k,$ then

\begin{eqnarray}\label{martingal1}
\nonumber
&&\mbox{with }
I(s) := N^2 (n_{k+1}(s) - n_k(s)) - N^2 (n_k(s) - n_{k-1}(s))
\\
\nonumber
&&\quad\quad\quad\quad
+\frac{N}{3 l} (n_{k+1}^2(s)+n_{k+1}(s) n_k(s)+n_k^2(s))
\\
\nonumber
&&\quad\quad\quad\quad
-\frac{N}{3 l} (n_k^2(s)+n_k(s) n_{k-1}(s)+n_{k-1}^2(s)),
\\
&&\nkN(t)-\int_0^t I(s) d s\mbox{ is a } {\cal
F}^N_t\mbox{-Martingale.}
\end{eqnarray}

Note that first, with a stopping time
$\tau_M$ such that
$\sup_{0\le t \le T} \sup_{k=1}^N n^N_k(t\wedge\tau_M)
1_{\{\tau_M>0\}} < M$,
we obtain that
$\nkN(t\wedge\tau_M)-\int_0^{t\wedge\tau_M} I(s) d s$ is a
${\cal F}^N_t$-Martingale for all $M>0.$ Equation
(\ref{martingal1}) will then follow from the proof of Lemma
\ref{S8b}.
Our approximating Markov jump process will be

\begin{equation}\label{process}
X^N(t,r) := X^{N, l}(t,r) := \frac{\nkN(t)}{l},\mbox{ }r\in
[\frac{k-1}{N},\frac{k}{N}),\mbox{ (and periodic extension).}
\end{equation}

Let $H^N$ be the $L^2(0,1)$-subspace of step functions on $[0,1)$
which are constant on the intervals $[\frac{k-1}{N},\frac{k}{N}).$
Define the orthogonal projection $P_N: L^2(0,1) \to H^N$ by

\begin{equation}\label{PN}
P_N f(r) = N \int_{\frac{k-1}{N}}^{\frac{k}{N}}f(x) d x
\mbox{ for }r\in[\frac{k-1}{N},\frac{k}{N}),
\end{equation}

and introduce the discrete derivatives

\begin{eqnarray}\label{discretederiv}
\nonumber \nabla_N^{\pm} f(r) &=& \pm N [ P_N f(r \pm N^{-1}) - P_N
f(r) ],\\
\Delta_N f(r) &=& \nabla_N^- \nabla_N^+ f(r) =
\nabla_N^+ \nabla_N^- f(r)\\
\nonumber&=& N^2[P_Nf(r+N^{-1})-2P_N f(r) + P_N f(r-N^{-1})].
\end{eqnarray}

From (\ref{process}) and (\ref{martingal1}) follows that

\begin{equation}\label{processint}
X^N(t) = X^N(0) + \int_0^t \Delta_N X^N(s) +
\nabla_N^+ F_N(X^N(s)) d s + Z^N(t)
\end{equation}

where

\begin{equation}\label{FN}
F_N: \left\{\begin{array}{l} H^N \mapsto H^N \\
X \to \frac13[
(X(\cdot))^2+X(\cdot)X(\cdot-N^{-1})+(X(\cdot-N^{-1}))^2]
\end{array}\right.
\end{equation}

and $Z^N(t)$ is an $H^N$-valued martingale for ${\cal F}^N_t.$
In mild form this becomes

\begin{equation}\label{processmild}
X^N(t) = e^{\Delta_N t} X^N(0) + \int_0^t e^{\Delta_N (t-s)}
\nabla_N^+ F_N (X^N(s)) d s + Y^N(t)
\end{equation}

where

\begin{equation}\label{convolution}
Y^N(t) = \int_0^t e^{\Delta_N (t-s)} d Z^N(s)
\end{equation}

(note that $Z^N$ is of bounded variation $P$-a.s. because it is
piecewise absolutely continuous). For technical reasons we assume that $N$ is
odd.

We obtain the following result (for the definition of the spaces
see Definition \ref{spaces}).
$\Skalprod{\cdot}{\cdot}$ denotes the dual pairing
and simultaneously the $L^2(0,1)$-scalar product.

\begin{Thm}\label{Thmstoch}
Let $X^N$ be the process defined by (\ref{process}) with
deterministic initial condition $0 \le X^N(0) \in H^N,$ such that for
arbitrary $\alpha\in (0,\frac12)$
$$
\Vert X^N(0) - \psi_0 \Vert_{H^\alpha(0,1)}
\stackrel{N \to \infty}{\longrightarrow} 0,
$$
where $0 \le \psi_0 \in H^\alpha(0,1)$ is the initial condition of
(\ref{stochburgeq}).
Moreover assume $l \ge q N$ for arbitrary $q>0.$

Then there exists a probability space
$(\tilde\Omega,\tilde{\cal F}, \tilde P),$ subsequences
$(N_k)_{k\in\N}$ and $(l_k)_{k\in\N},$ and $H^{N_k}$-valued processes
$\tilde X^{N_k},$ $\tilde Y^{N_k},$ and $\tilde Z^{N_k}$ on this
probability space. The common distribution of $\tilde X^{N_k},$
$\tilde Y^{N_k},$ and $\tilde Z^{N_k}$ equals the common distribution
of $X^{N_k},$ $Y^{N_k},$ and $Z^{N_k},$ for each $k \in \N.$
There exist processes
$\psi$ in $C(0,T,L^2(0,1)),$ $\tilde Y$ in
$C(0,T,H^{\alpha_1}(0,1)),$ and $M$ in
$C(0,T,H^{\alpha_2}(0,1))$ with $\alpha_1 < \frac12$ and
$\alpha_2 < - \frac12.$ $M$ is a martingale w.r.t.
$(\sigma( \psi(s),s \le t))_t.$ We obtain
$$
(\tilde X^{N_k}, \tilde Y^{N_k}, \tilde Z^{N_k})
\stackrel{k \to \infty}{\longrightarrow} (\psi, \tilde Y,
M)
$$
$\tilde P$-almost sure in
$L^2(0,T,L^2(0,1)) \times D(0,T,H^{\alpha_1}(0,1)) \times
D(0,T,H^{\alpha_2}(0,1)).$ The equation

\begin{equation}\label{stochburgmild}
\psi(t) = e^{t \Delta} \psi(0) +
\int_0^t e^{(t-s)\Delta} \nabla (\psi(s))^2 d s +
\tilde Y(t)
\end{equation}

is fulfilled $\tilde P$-a.s. in $C(0,T,L^2(0,1))$
where
$\tilde Y(t)= \int_0^t e^{(t-s)\Delta} d M(s)$
$\tilde P$-a.s. in $C(0,T,H^{\alpha_1}(0,1)).$
Here $e^{t\Delta}$ denotes the semigroup
defined by the Laplacian $\Delta$ with periodic boundary
conditions. The equation

\begin{equation}\label{stochburgschwach}
\Skalprod{\psi(t)}{\varphi}=\Skalprod{\psi(0)}{\varphi} +
\int_0^t \Skalprod{\Delta \psi(s) + \nabla (\psi(s))^2}{\varphi} d s
+ \Skalprod{M(t)}{\varphi}
\end{equation}

holds $\tilde P$-a.s. in $C(0,T,\R)$ where
$\varphi \in C^{\alpha_3}_{per}(0,1)$ with $\alpha_3>\frac52$ and
$\Skalprod{M(t)}{\varphi}
= \int_0^t \int_0^1 \sqrt{\gamma \psi(s,x)}
\varphi(x) d W(s,x)$ where $W$ is a certain space-time-white noise.
In this sense, $M$ can be represented as
$$M(t)= \int_0^t \sqrt{\gamma \psi(s)} d W(s).$$
\end{Thm}

\begin{Rem}
The proof of the theorem will be given in the next section in
Lemmata \ref{S16b}, \ref{S16c}, \ref{S16d}, \ref{S18}, and
\ref{S19a}. We have not tried to prove uniqueness of a solution of
Equation (\ref{stochburgeq}), see \cite{Bonnet}. We can generalize
Theorem \ref{Thmstoch} to random initial conditions. Note that the
$l$-dependence of the quantities in the theorem is suppressed in
the notation.
\end{Rem}

\section{Proofs}\label{zwei}

\begin{Defi}\label{eigen}

\noindent
(i) Eigenfunctions of $\Delta:$ Set $\varphi_0(r) := 1$ and
\begin{eqnarray}\label{phin}
\nonumber
&&
\varphi_n(r) := \sqrt{2} \sin(2\pi n r)\mbox{ for } n\in\N,
\\
&&
\varphi_n(r) := \sqrt{2} \cos(2\pi n r)\mbox{ for }
n\in\Z\setminus\N_0.
\end{eqnarray}
The eigenfunctions of $\Delta$ with periodic boundary conditions
on $(0,1)$ corresponding to the eigenvalues $\lambda_n = -4\pi^2n^2$
are given by the complete orthonormal system
$(\varphi_n)_{n\in\Z}\subset L^2(0,1).$

\bigskip\noindent
(ii) Eigenfunctions of $\Delta_N:$ Let
\begin{equation}\label{phinN}
\varphi_{n,N}(r) := \varphi_n(\frac{k-1}{N})\mbox{ for }
r \in [\frac{k-1}{N},\frac{k}{N}),
\end{equation}
where $k=1,\ldots,N, n=-\frac{N-1}{2},\ldots,\frac{N-1}{2}$ and
$N$ is assumed to be odd. According to \cite{Blount96},
$(\varphi_{n,N})_n$ form a complete orthonormal system in the space
$H^N\subset L^2(0,1)$ of piecewise constant functions (defined in
Section \ref{eins}). They are the eigenfunctions of $\Delta_N$
corresponding to the eigenvalues $\beta_{n,N} =
-2N^2(1-\cos(\frac{2\pi n}{N})).$ There are constants $0< c_1 <
c_2$ with
\begin{equation}\label{lb}
c_1 |\lambda_n| < |\beta_{n,N}| < c_2 |\lambda_n|
\end{equation}
for all $n=-\frac{N-1}{2},\ldots,\frac{N-1}{2}.$

\bigskip\noindent
(iii) Projection operators: $P_N$ is the
$L^2(0,1)$-orthogonal projection on $H^N$ and $P_n$ the
$L^2(0,1)$-orthogonal projection on
span$\{\varphi_k,k=-n,\ldots,n\}.$
\end{Defi}

\begin{Defi}\label{spaces}
We define the usual Sobolev spaces of order $\alpha\in\R$ with
periodic boundary conditions by
$$
H^\alpha(0,1) := \{ f = \sum_{n\in\Z} \alpha_n
\varphi_n, (\alpha_n)_{n\in\Z} \subset \R
\mbox{ with } \Vert f\Vert_{H^\alpha(0,1)} < \infty\}
$$
where $\Vert f\Vert_{H^\alpha(0,1)}^2 := \sum_{n\in\Z} \alpha_n^2
(1-\lambda_n)^\alpha.$ Similarly we set
$$
H^\alpha_N(0,1) := \{ f \in H^N:
\Vert f\Vert^2_{H^\alpha_N(0,1)} :=
\sum_{n=-\frac{N-1}{2}}^{\frac{N-1}{2}}
\Skalprod{f}{\varphi_{n,N}}^2 (1-\beta_{n,N})^\alpha < \infty\}.
$$
\end{Defi}

Set $\delta X(t) = X(t)-X(t-) = X(t) -
\lim_{s<t,s\to t} X(s).$ Then the following are ${\cal
F}_t^N$-martingales.
\begin{eqnarray}\label{martingalz}
\nonumber
&&
Z^N_D(t) := \sum_{s\le t} \delta X^N_D(s) - \int_0^t \Delta_N X^N(s)
+ \nabla_N^+ F_N(X^N(s)) d s,
\\
&&
Z^N_B(t) := \sum_{s\le t} \delta X^N_B(s),
\end{eqnarray}
where $\delta X^N_D \in H^N$ is a jump caused by diffusion and
$\delta X^N_B$ is a jump by birth or death. The proof is similar
to \cite{Blountdiss}. Moreover,
\begin{eqnarray}\label{martingalq}
\nonumber
&&\Skalprod{Z^N_D(t)}{f}^2 - \frac{1}{N l}\int_0^t
\Skalprod{X^N(s)}{(\nabla_N^+f)^2}
+\Skalprod{X^N(s)+\frac{1}{N}F_N(X^N(s))}{(\nabla_N^-f)^2} d s,
\\
&&\Skalprod{Z^N_B(t)}{f}^2-\gamma  \int_0^t
\Skalprod{X^N(s)}{f^2} d s
\end{eqnarray}
are ${\cal F}_t^N$-martingales, $f \in H^N.$

\begin{Lemma}\label{S8b}
Let the conditions of Theorem \ref{Thmstoch} be fulfilled. Then
with $Y^N_B(t) = \int_0^t e^{\Delta_N (t-s)} d Z^N_B(s),$
$$
\sup_{N} P(\Vert Y_B^N \Vert_{L^\infty(0,T,H_N^{\alpha_1})}
\ge \tilde R) \stackrel{ \tilde R \to \infty}{\longrightarrow} 0,
$$
for $\alpha_1<\frac12.$
\end{Lemma}

\begin{Bew} The proof follows \cite{Blount96}, Lemma
3.2. We therefore only give a brief sketch of the idea. Let
$$
R(t) = \int_0^t e^{(t-s)\Delta_N} d Z_B^N(s \wedge \tau)
$$
where
$\tau = \tau_N
:= \inf\{t \in [0,T]: \Skalprod{X^N(t)}{1} \ge \rho\}$.
Since
$P(\Vert Y_B^N \Vert_{L^\infty(0,T,H_N^{\alpha_1})}
\ge \tilde R)
\le
P(\Vert R \Vert_{L^\infty(0,T,H_N^{\alpha_1})}
\ge \tilde R)
+
P(\tau_N < T)$
we have to show $\sup_N P(\tau_N < T)
\stackrel{ \rho \to \infty}{\longrightarrow} 0$ and
for fixed $\rho>0,$
$\sup_N P(\Vert R \Vert_{L^\infty(0,T,H_N^{\alpha_1})}
\ge \tilde R)
\stackrel{ \tilde R \to \infty}{\longrightarrow} 0$.
Let now $\rho>0$ be fixed and
define for $m \not= 0$ and $u \in [0,t]$
$$
M(u) = |m| \int_0^u e^{\beta_{m,N}(t-s)}
d \Skalprod{Z_B^N(s \wedge \tau)}{\varphi_{m,N}}.
$$
This is a mean-zero-martingale with
$M(t) = |m| \Skalprod{R(t)}{\varphi_{m,N}}$ and
$|\delta M(u)| \le 1.$ The predictable quadratic variation process $<<M>>$
fulfills $<<M>>(u) \le c\gamma \rho,$ see (\ref{martingalq}).
Lemma 4.4 of
\cite{Blount92} yields $E[\exp(M(t))] \le \exp(\frac32 c \gamma
\rho)$ whence
$$
P(m^{2\alpha_1} \Skalprod{R(t)}{\varphi_{m,N}}^2 \ge m^{-2r})
\le c(\gamma \rho) \exp(-|m|^{1-r-\alpha_1})
$$
Because for $\alpha_1 < \frac12$ there exists $r > \frac12$ with
$\alpha_1+r < 1$ such that $\sum_{m \in \Z\setminus\{0\}}
|m|^{-2r} < \infty$ and
$\sum_{m \in \Z\setminus\{0\}} m^2 \exp(-|m|^{1-r-\alpha_1})
< \infty,$ we obtain as in \cite{Blount96}
\begin{eqnarray*}
&&
\sup_N P(\sup_{t \le T} \Vert R(t) \Vert_{H^{\alpha_1}_N} \ge
\tilde R)
\le c(\gamma,\rho,T,\alpha_1)
\sum_{m \in \Z\setminus\{0\}}
m^2 \exp(-c(T)|m|^{1-r-\alpha_1}\tilde R)
\\
&&
+
\sup_N P(\sup_{t \le T} \Skalprod {R(t)}{1}
\ge \frac{\tilde R}{c})
\stackrel{ \tilde R \to \infty}{\longrightarrow} 0
\end{eqnarray*}
where the last term is treated similarly to the others.
For the assertion of the Lemma it therefore remains to show
$\sup_N P(\tau_N < T)
\stackrel{ \rho \to \infty}{\longrightarrow} 0$ which in turn
follows from $E[\sup_{t \le T} \Skalprod{X^N(t)}{1}] \le c$
uniformly in $N:$ From
$\Skalprod{\Delta_N X^N(t)+\nabla_N^+ F_N(X^N(t))}{1} = 0$
we conclude
$\Skalprod{X^N(t)}{1}
= \Skalprod{X^N(0)}{1} + \Skalprod{Z^N(t)}{1}$ and by the
Jensen and maximal inequality and (\ref{martingalz})
\begin{eqnarray*}
&&
E[\sup_{t\le T} \Skalprod{X^N(t \wedge \tau)}{1}]
= E[\sup_{t\le T} \Skalprod{Z^N(t \wedge \tau)}{1}]
+ \Skalprod{X^N(0)}{1}
\\
&&
\le
2 \sup_{t\le T}
\sqrt{ E[ \Skalprod{Z^N(t \wedge \tau)}{1}^2]}
+ \Skalprod{X^N(0)}{1}
\\
&&
\le
4 \sup_{t\le T}
\sqrt{ E[ \Skalprod{Z^N_D(t \wedge \tau)}{1}^2]
+ E[ \Skalprod{Z^N_B(t \wedge \tau)}{1}^2]}
+ \Skalprod{X^N(0)}{1} \le
\end{eqnarray*}
since $\Skalprod{Z^N_D(t \wedge \tau)}{1}=0$ a.s, we
continue using Equation (\ref{martingalq}) and
$E[\Skalprod{Z^N(t \wedge \tau)}{1}]=0$
$$
\le 4 \sup_{t\le T} \sqrt{
\gamma \int_0^{t \wedge \tau} \Skalprod{X^N(0)}{1}
+ E[\Skalprod{Z^N(s \wedge \tau)}{1}] d s } + \Skalprod{X^N(0)}{1}
\le c
$$
uniformly in $N.$
\end{Bew}

\begin{Lemma}\label{S10}
Let the conditions of Theorem \ref{Thmstoch} be fulfilled. Then
with $Y^N_D(t) = \int_0^t e^{\Delta_N (t-s)} d Z^N_D(s),$
$$
\sup_{N} P(\Vert Y_D^N \Vert_{L^\infty(0,T,H_N^{\alpha_1})}
\ge \tilde R) \stackrel{ \tilde R \to \infty}{\longrightarrow} 0,
$$
for $\alpha_1<\frac12.$
\end{Lemma}

\begin{Bew}
We proceed as in
the proof of Lemma \ref{S8b}. Due to
(\ref{martingalq}) we obtain with the notation of this proof:
\begin{eqnarray*}
&&
<<M>>(u) = \frac{m^2}{N l} \int_0^{u\wedge \tau}
\exp( 2\beta_{m,N}(t-s))
\\
&&
\quad\quad
\times \Big(
\Skalprod{X^N(s)}{(\nabla_N^+ \varphi_{m,N})^2}
+\Skalprod{X^N(s)+\frac{1}{N}F_N(X(s))}
{(\nabla_N^- \varphi_{m,N})^2} \Big) d s \le
\end{eqnarray*}
by $\Vert\nabla_N^+ \varphi_{m,N}\Vert_{L^\infty} \le c m$
$$
\le c \frac{m^4}{N l} \int_0^{u\wedge \tau}
\exp( 2\beta_{m,N}(t-s)) \Big( \Vert X^N(s)\Vert_{L^1}
+ \frac{1}{N}\Vert X^N(s) \Vert^2_{L^2} \Big) d s \le
$$
by $\Vert X^N(s) \Vert^2_{L^2} \le N \Vert X^N(s)\Vert_{L^1}^2$
$$
\le c \frac{m^2}{N l}(\rho+\rho^2) \le c(\rho).
$$
We can continue similarly to the proof of Lemma \ref{S8b}.
\end{Bew}

\begin{Lemma}\label{S14}
Under the assumptions of Theorem \ref{Thmstoch}, the family of the
probability distributions of $Y^N_B$ is tight on
$D(0,T,H^{\alpha_1}(0,1)).$
\end{Lemma}

\begin{Bew}
We again follow the proof of \cite{Blount96}, Lemma 3.3 and first
show that the distributions of $Z^N_B$ are tight on
$D(0,T,H^{\alpha_2}(0,1)).$ We verify Condition (a) in Theorem 37.2
in \cite{Ethier} and (8.33) and (8.29) ibid. Let $\alpha_2 <
\tilde \alpha < -\frac12$ and $\Gamma_{\eta}:= B_R(0) \subset
H^{\tilde \alpha}(0,1) \subset \subset H^{\alpha_2}(0,1).$ Then by
(\ref{martingalq})
\begin{eqnarray*}
&&
P(Z_B^N(t) \in \Gamma_{\eta})
\ge 1-P(\Vert Z_B^N(t)\Vert^2_{H^{\tilde \alpha}_N} > R^2)
\\
&&
\ge 1- \frac{\gamma}{R^2} \sum_{m=-\frac{N-1}{2}}^{\frac{N-1}{2}}
\int_0^t E[\Skalprod{X^N(s)}{\varphi^2_{m,N}}]
(1-\beta_{m,N})^{\tilde\alpha} d s \ge 1-\eta
\end{eqnarray*}
for sufficiently large $R,$ see the proof of Lemma \ref{S8b}.
Analogously, for $0 \le t \le T$ and $0 \le u \le \tilde \delta$
for some $\tilde \delta \in (0,1),$ $E[\Vert Z^N_B(t+u) - Z^N_B(t)
\Vert^2_{H^{\tilde\alpha}(0,1)} | {\cal F}^N_t] \le c \tilde\delta.$

This shows the existence of
$K_\eta \subset \subset D(0,T,H^{\alpha_2}(0,1))$
with $P(Z^N_B \in K_\eta) \ge 1-\eta.$ By Lemma \ref{d2p},
$K_\eta \subset \subset L^p(0,T,H^{\alpha_2}(0,1))$ for all $p>1.$
According to Lemma \ref{S12},
\begin{eqnarray*}
&&
\{ \int_0^t e^{(t-s)\Delta_N}\Delta_N Z^N_B(s) d s \in
C(0,T,H^{\alpha_2-\epsilon}(0,1)) | Z^N_B \in K_\eta\}
\\
&&
\subset\subset C(0,T,H^{\alpha_2-\epsilon}(0,1))
\subset D(0,T,H^{\alpha_2-\epsilon}(0,1))
\end{eqnarray*}
for all $\epsilon > 0.$ Because $P_n$ is continuous from
$D(0,T,H^{\alpha_2-\epsilon}(0,1))$ into\\
$D(0,T,H^{\alpha_1}(0,1))$, and because of
\begin{equation}\label{partint}
Y_B^N(t)=Z_B^N(t) +
\int_0^t e^{(t-s)\Delta_N}\Delta_N Z^N_B(s) d s,
\end{equation}
see \cite{Blount96}, the distributions of $P_n Y_B^N$ are tight
on
$D(0,T,H^{\alpha_1}(0,1))$ for fixed $n\in\N.$ The assertion of
the Lemma then
follows from Problem 18, Chapter 3, in \cite{Ethier} and the fact
that for all $\epsilon > 0$ there exist a $n\in\N$ such that
$P(\Vert (I-P_n) Y^N_B \Vert_{L^\infty(0,T,H^{\alpha_1}(0,1))}
\ge \epsilon) \le \epsilon$ uniformly in $N.$ This can
easily be deduced from Lemma \ref{S8b} and Lemma \ref{HaN}.
\end{Bew}

\begin{Lemma}\label{S16a}
Let the requirements of Theorem \ref{Thmstoch} be true. Then $Y^N_D
\to 0$ for $N\to\infty$ in $L^\infty(0,T,L^2(0,1))$ in
probability.
\end{Lemma}

\begin{Bew}
According to Lemma \ref{S10} it suffices to show
$P_n Y^N_D( \cdot \wedge \tau ) \to 0$ for $N\to\infty$ in the
$L^\infty(0,T,L^2(0,1))$-norm in probability for any fixed $n$, where
$\tau$ is the stopping time from the proof of Lemma \ref{S8b}. Proving
this can be done as in \cite{Blount91}.
\end{Bew}

Let
\begin{equation}\label{un}
u^N=X^N-Y^N.
\end{equation}
By (\ref{processmild}), $u^N:[0,T]\to H^N$ is
continuous. Moreover, between two jumps, we obtain
\begin{equation}\label{un2}
\dnachd{t} u^N(t) = \Delta_N u^N(t) + \nabla_N^+ F_N(u^N+Y^N)(t).
\end{equation}

\begin{Lemma}\label{S6}
With constants $c$ independent of $N,$ $u^N,$ and $Y^N$ we obtain
the following a-priori-estimates:
\begin{eqnarray*}
&&
\Vert u^N(t)\Vert^2_{L^2(0,1)}
\le
\exp\Big(c\int_0^t(1+\Vert Y^N(s) \Vert^{8/3}_{L^4(0,1)})d s\Big)
\\
&&
\quad\quad\quad\quad\quad\quad\quad
\times
\Big(\Vert X^N(0)\Vert^2_{L^2(0,1)}
+ c \int_0^t\Vert Y^N(s)\Vert^{4}_{L^4(0,1)} d s\Big) =:f(t)
\end{eqnarray*}
and
\begin{eqnarray*}
&&
\int_0^T \Vert \nabla^-_N u^N(t)\Vert^2_{L^2(0,1)} d t
\\
&&
\le \int_0^T \Big( c f(t) (1+\Vert Y^N(t)\Vert^{8/3}_{L^4(0,1)})
+c \Vert Y^N(t)\Vert^{4}_{L^4(0,1)} \Big) d t
+ \Vert X^N(0)\Vert^2_{L^2(0,1)}.
\end{eqnarray*}
\end{Lemma}

\begin{Bew}
We apply a well-known procedure, see e.g. \cite{Prato}.
Both estimates follow from
\begin{equation}\label{diffapriori}
\dnachd{t} \Vert u^N\Vert^2_{L^2}
+ \Vert \nabla^-_N u^N \Vert^2_{L^2}
\le c \Vert u^N\Vert^2_{L^2}
\Big(1+\Vert Y^N\Vert^{8/3}_{L^4}\Big) + c \Vert Y^N\Vert^4_{L^4}.
\end{equation}
by an application of the Gronwall lemma. To obtain this, we
multiply (\ref{un2}) with $u^N$ and integrate over the
spatial variable,
$$
\frac12 \dnachd{t} \Skalprod{u^N}{u^N}
- \Skalprod{\Delta_N u^N}{u^N}
= -\Skalprod{\nabla_N^+ F_N(u^N+Y^N)}{u^N}.
$$
By partial integration
$\Skalprod{\nabla_N^+ f}{g} = - \Skalprod{f}{\nabla_N^- g}$
for $f, g \in H^N$ we obtain
\begin{equation}\label{UN2}
\frac12 \dnachd{t} \Vert u^N \Vert_{L^2(0,1)}^2
+ \Vert \nabla^-_N u^N \Vert_{L^2(0,1)}^2 = S_1 + S_2 + S_3.
\end{equation}
with
\begin{eqnarray*}
&&
S_1 = -\frac13\Skalprod
{(u^N(\cdot))^2+u^N(\cdot)u^N(\cdot-\frac{1}{N})
+(u^N(\cdot-\frac{1}{N}))^2}
{\nabla^-_N u^N},
\\
&&
S_2 = -\frac13\Skalprod
{2u^N(\cdot)Y^N(\cdot)
+2u^N(\cdot-\frac{1}{N})Y^N(\cdot-\frac{1}{N})
\\
&&
\quad\quad\quad+u^N(\cdot)Y^N(\cdot-\frac{1}{N})
+u^N(\cdot-\frac{1}{N})Y^N(\cdot)}
{\nabla^-_N u^N},
\\
&&
S_3 = -\frac13\Skalprod
{(Y^N(\cdot))^2+Y^N(\cdot)Y^N(\cdot-\frac{1}{N})
+(Y^N(\cdot-\frac{1}{N}))^2}
{\nabla^-_N u^N}.
\end{eqnarray*}
We now treat $S_1,$ $S_2,$ and $S_3.$

$S_1=0,$ because
$$
S_1=-\frac{N}{3} \sum_{k=1}^{N}
(u^N(\frac{k}{N}))^3-(u^N(\frac{k-1}{N}))^3=0
$$
due to periodic boundary conditions.
This is the discrete equivalent to the standard trick
$\int_0^1 u^2 \dnachd{x} u =  \frac13 \int_0^1 \dnachd{x} u^3 =
0$ used when treating the Burgers equation with periodic boundary
conditions. Note that for this result we impose the jump rate in
(\ref{jumprates}). Otherwise we could have taken
$N^2 n_k + \frac{N}{l} n_k^2$ as first rate in (\ref{jumprates})
which would have entailed $F_N(X) = X^2$ in (\ref{FN}) and
therefore a ``usual'' deterministic Burgers equation in
(\ref{processmild}) and (\ref{un2}).

$S_2+S_3$ can be bounded by
$c\Vert u^N\Vert_{L^4} \Vert Y^N\Vert_{L^4}
\Vert\nabla_N^-u^N\Vert_{L^2}+\Vert Y^N\Vert_{L^4}^2
\Vert\nabla_N^-u^N\Vert_{L^2}.$
Since
$\Vert u^N\Vert_{L^4} \le c \Vert u^N\Vert_{H^{\frac14}_N}$
due to a standard Sobolev imbedding (\cite{Adams}, Theorem 7.57)
and Lemma \ref{HaN}, and
due to Lemma \ref{S7},
$$
\Vert u^N\Vert_{H^{\frac14}_N}
\le c \Vert u^N\Vert^{\frac34}_{L^2}
\Vert u^N\Vert^{\frac14}_{H^1_N}
\le c \Vert u^N\Vert_{L^2} + c \Vert u^N\Vert^{\frac34}_{L^2}
\Vert\nabla_N^-u^N\Vert^{\frac14}_{L^2},
$$
the crucial term in the bound of $S_2+S_3$ is
$$c\Vert u^N\Vert^{\frac34}_{L^2}
\Vert\nabla_N^-u^N\Vert^{\frac54}_{L^2} \Vert Y^N\Vert_{L^4}
\le \frac16 \Vert\nabla_N^-u^N\Vert^2_{L^2}
+ c \Vert u^N\Vert^2_{L^2} \Vert Y^N\Vert^{\frac83}_{L^4}.
$$ This yields (\ref{diffapriori}).
\end{Bew}

\begin{Lemma}\label{S13}
Let the conditions of Theorem \ref{Thmstoch} be fulfilled. Then
the family of the probability distributions of $u^N$ is tight on
$L^2(0,T,L^2(0,1)).$
\end{Lemma}

\begin{Bew}
By the computation in the proof of Lemma \ref{S7} we have for all
$\beta\in\R$ and $v^N\in H^N,$
$\Vert \nabla_N^+ v^N\Vert^2_{H^\beta_N} \le c
\Vert v^N\Vert^2_{H^{\beta+1}_N}.$
A rough estimate now gives
\begin{eqnarray*}
&&
\Vert F_N(u^N+Y^N)\Vert^2_{H^{\beta+1}_N}
\le c
\sum_{m=1}^\frac{N-1}{2}
\Skalprod{F_N(u^N+Y^N)}{1}^2(1-\beta_{m,N})^{1+\beta}
\\
&&
\le c \Vert u^N+Y^N \Vert^4_{L^2}
\sum_{m=1}^\frac{N-1}{2} (1-\beta_{m,N})^{1+\beta}
\le c \Vert u^N+Y^N \Vert^4_{L^2}
\end{eqnarray*}
for $\beta < - \frac32.$ Therefore
\begin{eqnarray*}
&&
\sup_N P\Big( \int_0^T
\Vert \nabla_N^+F_N(u^N+Y^N)\Vert^2_{H^\beta_N} d t \ge R \Big)
\\
&&
\le
\sup_N P\Big( \int_0^T \Vert u^N+Y^N \Vert^4_{L^2}
\ge\frac{R}{c}\Big) \stackrel{R\to \infty}{\longrightarrow}0
\end{eqnarray*}
according to Lemmata \ref{S6}, \ref{S8b} with Lemma \ref{HaN}, and
Lemma \ref{S10}. Let $(R_N h)(t) = \int_0^t e^{(t-s)\Delta_N} h(s)
d s$
for $h\in H^N$ and
\begin{eqnarray*}
\Xi(R) &=& \{ u^N \in C(0,T,H^\beta) \cap L^2(0,T,H^1_N):
\\
&&
\Vert \nabla_N^+F_N(u^N+Y^N)\Vert_{L^2(0,T,H^\beta_N)} \le R,
\Vert \nabla_N^+ u^N\Vert_{L^2(0,T,L^2)}\le R,
\\
&&
u^N(t) = e^{t\Delta_N} X^N(0)
+ R_N(\nabla_N^+F_N(u^N+Y^N))(t)\}.
\end{eqnarray*}
Then by the equivalence of the norms in $H^\beta$ and $H^\beta_N$
for $\beta\le 0,$
\begin{eqnarray*}
&&
P(\Xi(R))
\\
&&
\ge
1-P(\Vert \nabla_N^+F_N(u^N+Y^N)\Vert_{L^2(0,T,H^\beta)} \ge R)
-P(\Vert \nabla_N^+ u^N\Vert_{L^2(0,T,L^2)}\ge R)
\\
&&
\ge 1-\epsilon
\end{eqnarray*}
for $R=R(\epsilon)$ according to Lemmata \ref{S6}, \ref{S8b} with
Lemma \ref{HaN}, and Lemma \ref{S10}. Moreover,
$\Xi(R)$ is compact in $L^2(0,T,L^2)$ according to Lemmata
\ref{S12} and \ref{S11}.
\end{Bew}

\begin{Lemma}\label{S16b}
Let the conditions of Theorem \ref{Thmstoch} be fulfilled. Then
there exist subsequences $(N_k)_{k\in\N}$ and $(l_k)_{k\in\N}$ and
a probability measure $\mu$ such that in distribution on
$L^2(0,T,L^2(0,1)) \times D(0,T,H^{\alpha_1}(0,1)) \times
D(0,T,H^{\alpha_2}(0,1)) \times D(0,T,H^{\alpha_1}(0,1)),$
$$
(u^{N_k}, Y_B^{N_k}, Z_B^{N_k}, Y_D^{N_k})
\stackrel{k \to \infty}{\longrightarrow} \mu
$$
for $\alpha_1 < \frac12$ and $\alpha_2 < - \frac12.$
\end{Lemma}

\begin{Bew} See Lemma \ref{S16a}, Lemma \ref{S14} and its proof,
and Lemma \ref{S13}.  The tightness of the family of the
probability distributions of $Y^N_D$ is shown as in
Lemma \ref{S14} using an estimate
$E[ \Vert X^N(t)\Vert^2_{L^2(0,1)}] \le c_T N$ derived similarly to
the proof of Lemma \ref{S8b}.
The conclusion follows from the theorem of
Prokhorov, e.g. \cite{Ethier}, Chapter 3.
\end{Bew}

\begin{Lemma}\label{S16c}
Let the requirements of Theorem \ref{Thmstoch} be fulfilled. There
exists a probability space
$(\tilde\Omega,\tilde{\cal F}, \tilde P)$ and processes
$\tilde u$ in $L^2(0,T,L^2(0,1)),$ $\tilde Y_B$ in
$D(0,T,H^{\alpha_1}(0,1)),$ and $M$ in
$D(0,T,H^{\alpha_2}(0,1))$ with $\alpha_1 < \frac12$ and
$\alpha_2 < - \frac12.$
$(\tilde u,\tilde Y_B,M,0)$ has the law defined in Lemma \ref{S16b}.
There exist processes $\tilde u^{N_k},$
$\tilde Y_B^{N_k},$ $\tilde Y_D^{N_k},$ and $\tilde Z_B^{N_k}$ on
this probability space such that the common distribution of
$\tilde u^{N_k},$ $\tilde Y_B^{N_k},$ $\tilde Y_D^{N_k},$ and
$\tilde Z_B^{N_k}$ equals the common distribution of $u^{N_k},$
$Y_B^{N_k},$ $Y_D^{N_k},$ and $Z_B^{N_k},$ for each $k \in \N.$
Moreover, in $L^2(0,T,L^2(0,1)) \times D(0,T,H^{\alpha_1}(0,1))
\times D(0,T,H^{\alpha_2}(0,1)) \times D(0,T,H^{\alpha_1}(0,1)),$
$$
(\tilde u^{N_k}, \tilde Y_B^{N_k}, \tilde Z_B^{N_k}, \tilde
Y_D^{N_k}) \stackrel{k \to \infty}{\longrightarrow} (\tilde u,
\tilde Y_B, M, 0),
$$
$ \tilde P$-almost surely.
\end{Lemma}

\begin{Bew}
This lemma follows from Lemma \ref{S16b}, Lemma \ref{S16a} and the
theorem of Skorohod, e.g. \cite{Ethier}, Chapter 3.
\end{Bew}

\begin{Lemma}\label{S16d}
Let the assumptions of Lemma \ref{S16c} be fulfilled. Then
\begin{equation}\label{darst1}
\tilde u(t) = e^{t \Delta} \psi_0 +
\int_0^t e^{(t-s)\Delta}\nabla (\tilde u + \tilde Y_B)^2 (s) d s
\end{equation}
holds in $L^2(0,T,L^2(0,1)),$ $\tilde P$-a.s.
\end{Lemma}

\begin{Bew} For simplicity we denote the subsequence
$(N_k)_{k\in\N}$ in the Lemmata \ref{S16b} and \ref{S16c} by
$N,$ suppress the tilde, and replace $L^2(0,1)$ by $L^2,$ e.g.,
in this proof. We give the proof in several steps.

\bigskip\noindent
(i) $\Vert F_N(u^N)\Vert_{L^1} \le \Vert u^N\Vert^2_{L^2},$ see
(\ref{FN}).

\bigskip\noindent
(ii) $\Vert F_N(u^N) - u^2 \Vert_{L^1}
\le \Vert u^N - u\Vert_{L^2}
(\Vert u^N \Vert_{L^2}+\Vert u \Vert_{L^2})$ since
\begin{eqnarray*}
&&
\int_0^1 \frac13
\Big((u^N(x))^2+u^N(x)u^N(x-N^{-1})+(u^N(x-N^{-1}))^2\Big) -
\Big(u(x)\Big)^2 d x
\\
&\le&
 \frac13 \Big(
\Vert u^N+u\Vert_{L^2} \Vert u^N-u\Vert_{L^2}
+ \Vert u^N+u\Vert_{L^2}  \Vert u^N-u\Vert_{L^2}
\\
&&
+ \Vert u^N\Vert_{L^2} \Vert u^N(\cdot)-u^N(\cdot - N^{-1})\Vert_{L^2}
+ \Vert u^N+u\Vert_{L^2} \Vert u^N-u\Vert_{L^2}\Big)
\\
&\le&
\Vert u^N+u\Vert_{L^2} \Vert u^N-u\Vert_{L^2}
+\frac{1}{N} \Vert u^N\Vert_{L^2}  \Vert \nabla_N^- u^N\Vert_{L^2}.
\end{eqnarray*}

\bigskip\noindent
(iii) $\Vert \nabla \varphi - \nabla_N^- \varphi \Vert_{L^\infty}
\le \frac{1}{N} \Vert \varphi \Vert_{H^{a_1}}$ for $a_1>\frac52$
because:
\begin{eqnarray*}
&&
\Vert \nabla \varphi - \nabla_N^- \varphi \Vert_{L^\infty}
\\
&&
=\sup_{k=1,\ldots,N}\sup_{x\in [\frac{k-1}{N},\frac{k}{N})}
| \varphi'(x) - N
\int_{\frac{k-1}{N}}^{\frac{k}{N}}
\frac{\varphi(y)-\varphi(y-N^{-1})}{N^{-1}} d y |
\\
&&
=\sup_{k=1,\ldots,N}\sup_{x\in [\frac{k-1}{N},\frac{k}{N})}
| \varphi'(x) - \varphi'(\xi_{k,N}) |
\le \frac{c}{N} \Vert \varphi'' \Vert_{L^\infty}
\end{eqnarray*}
with some $\xi_{k,N}\in [\frac{k-1}{N},\frac{k}{N}].$

\bigskip\noindent
(iv) If $u^N \to u$ in $L^2(0,T,L^2)$ then $\nabla_N^+ F_N(u^N)
\to \nabla u^2$ in $L^1(0,T,H^{-a_1})$ for $a_1 > \frac52.$ For
\begin{eqnarray*}
&&
\Vert \nabla_N^+ F_N(u^N) - \nabla u^2\Vert_{H^{-a_1}}
=c\sup_{\Vert\varphi\Vert_{H^{a_1}}=1}
\Skalprod{\nabla_N^+ F_N(u^N) - \nabla u^2}{\varphi}
\\
&&
\le c \sup_{\Vert\varphi\Vert_{H^{a_1}}=1}
\Skalprod{F_N(u^N)}{\nabla \varphi - \nabla_N^- \varphi}
+ c \sup_{\Vert\varphi\Vert_{H^{a_1}}=1}
\Skalprod{u^2-F_N(u^N)}{\nabla\varphi}
\\
&&
\le \frac{c}{N} \Vert u^N \Vert^2_{L^2} + c
\Vert u^N- u \Vert_{L^2}
\Big(\Vert u^N\Vert_{L^2} + \Vert u\Vert_{L^2} \Big)
+\frac{c}{N} \Vert u^N\Vert_{L^2}  \Vert \nabla_N^- u^N\Vert_{L^2}
\end{eqnarray*}
according to steps (i), (ii), and (iii).

\bigskip\noindent
(v) For fixed $n \in N$,
$e^{\beta_{n,N}t} - e^{\lambda_n t} \to 0$ for $N \to \infty$
uniformly in $t \le T,$ because
$| e^{\beta_{n,N}t} - e^{\lambda_n t}|
\le t |\beta_{n,N} - \lambda_n|$
and $|\beta_{n,N} - \lambda_n|
= 4 \pi n^2 |2 \frac{1-\cos x}{x^2}-1| \to 0$ for $x \to 0$ where
$x=\frac{2\pi n}{N}.$

\bigskip\noindent
(vi) If $h^N$ is bounded in $L^1(0,T,H^{-a_1})$ then for
$a_2>a_1,$
$$
\Vert e^{\Delta_N t} h^N - e^{\Delta t} h^N
\Vert_{L^1(0,T,H^{-a_2})} \to 0 \mbox{ for } N \to \infty.
$$
Applying Lemma \ref{proj} we obtain
\begin{eqnarray*}
&&
\Vert e^{\Delta_N t} h^N - e^{\Delta t} h^N
\Vert^2_{H^{-a_2}}
\\
&&
=
\Vert \sum_{k\in\Z} \Big( \sum_{n=-\frac{N-1}{2}}^{\frac{N-1}{2}}
\Skalprod{h^N}{\varphi_{n,N}}\Skalprod{\varphi_{n,N}}{\varphi_k}
(e^{\beta_{n,N} t} - e^{\lambda_k t}) \Big) \varphi_k
\Vert^2_{H^{-a_2}}
\\
&&
\le
\sum_{n=-\frac{N-1}{2}}^{\frac{N-1}{2}} \sum_{l\in\Z}
(1-\lambda_{n+lN})^{-a_2}
\Skalprod {h^N}{\varphi_{n,N}}^2
(a^2_{n+lN}+b^2_{n+lN})
(e^{\beta_{n,N} t} - e^{\lambda_{n+lN} t})^2
\end{eqnarray*}
by a consideration similar to the proof of Lemma \ref{HaN},
distinguishing the cases $l=0,1,\ge 2$
\begin{eqnarray*}
&&
\le c \sum_{n=1}^{\frac{N-1}{2}}
\Big(
\Skalprod{h^N}{\varphi_{n,N}}^2 + \Skalprod{h^N}{\varphi_{-n,N}}^2
\Big)
\Big(
(1-\lambda_n)^{-a_2} |e^{\beta_{n,N}t} - e^{\lambda_n t}|^2\\
&&
\quad\quad\quad\quad
+ N^{2(a_1-a_2)} (1-\lambda_n)^{-a_1}
+ n^{-2 a_2} \frac{N}{n} (1-\cos(\frac{2\pi n}{N})
\Big)
\\
&&
\le c
\Big(
\max_{|n| \le k_0} |e^{\beta_{n,N}t} - e^{\lambda_n t}|^2
+ (1-\lambda_{k_0})^{a_1-a_2} + N^{2(a_1-a_2)}
\\
&&
\quad\quad\quad\quad
+ \max_{|n| \le k_0} |\frac{N}{n}(1-\cos(\frac{2\pi n}{N}))|
+ k_0^{2(a_1-a_2)}
\Big) \Vert h^N \Vert^2_{H^{-a_1}}.
\end{eqnarray*}
The claim follows from the boundedness of $h^N$ and (v).

\bigskip\noindent
(vii) For $X_0^N \stackrel{N \to \infty}{\longrightarrow} \psi_0$
in $H^\alpha,$ $\alpha>0,$ we obtain
$e^{\Delta_N t} X_0^N \to e^{\Delta t} \psi_0$ for $N\to\infty$ in
$L^\infty(0,T,L^2).$

\bigskip\noindent
(viii) For $u^N\stackrel{N \to \infty}{\longrightarrow}u,$ in
$L^2(0,T,L^2),$
$$
\int_0^t e^{(t-s)\Delta_N}\nabla_N^+ F_N(u^N)(s) d s
\stackrel{N \to \infty}{\longrightarrow}
\int_0^t e^{(t-s)\Delta}\nabla u^2 (s) d s
$$
in $L^1(0,T,H^{-a_2})$ because
\begin{eqnarray*}
&&
\int_0^T \Vert
\int_0^t e^{(t-s)\Delta_N}\nabla_N^+ F_N(u^N)(s) d s
- \int_0^t e^{(t-s)\Delta}\nabla u^2 (s) d s
\Vert_{H^{-a_2}} d t
\\
&&
\le \int_0^T \int_0^t
\Vert \Big( e^{(t-s)\Delta_N}-e^{(t-s)\Delta}\Big)
\nabla_N^+ F_N(u^N)(s) \Vert_{H^{-a_2}}
\\
&&
\quad\quad\quad\quad
+
\Vert e^{(t-s)\Delta}
\Big( \nabla_N^+ F_N(u^N)(s) - \nabla u^2 (s) \Big)
\Vert_{H^{-a_2}}
d s d t
\end{eqnarray*}
The first summand in the integral tends to $0$ because of (vi) and
(iv), the second due to (iv). See Lemma \ref{S6} and Lemma \ref{S16c}.

\bigskip\noindent
(ix) $u^N+Y_B^N+Y_D^N \stackrel{N \to \infty}{\longrightarrow}
u+Y_B,$ $\tilde P$-a.s. in $L^2(0,T,L^2)$ by Lemmata
\ref{S16c} and \ref{d2p}.
\end{Bew}

\begin{Lemma}\label{S18}
Under the requirements of Theorem \ref{Thmstoch}, $M$ and \\
$\Skalprod{M}{f}^2 - \gamma \int_0^t \Skalprod{\psi(s)}{f^2} d s$
are $(\sigma(\psi(s), s\le t))_t$-martingales, for all $f \in
L^\infty(0,1),$ where $\psi:=\tilde u + \tilde Y_B,$ see Lemma
\ref{S16c}. Moreover, $M \in C(0,T,H^{\alpha_2}(0,1)),$ $\tilde
P$-a.s. The quadratic variation process of $\Skalprod{M(t)}{f}$ is
given by $\gamma \int_0^t \Skalprod{\psi(s)}{f^2} d s,$ for all $f
\in L^\infty(0,1).$
\end{Lemma}

\begin{Bew}
We follow the proof of Lemma 3.6 of \cite{Blount96}. According to
the proof of Lemma \ref{S8b} we have, for $m \in \Z,$ that
$$
\tilde E[ \Skalprod{\tilde Z_B^{N_k}(t)}{\varphi_m}^2]
\le 2 \gamma \tilde E \int_0^t \Skalprod{\tilde X^{N_k}(s)}{1}
d s \le c
$$
uniformly in $k \in \N,$ where
$\tilde X^{N_k} := \tilde u^{N_k} + \tilde Y_B^{N_k}.$ From
\cite{Ethier}, Chapter 7, Problem 7, we infer that
$\Skalprod{M}{\varphi_m}$ and then $M$ are martingales w.r.t.
the above filtration.

By the Burkholder inequality and the proof
of Lemma \ref{S8b},
$$
\tilde E[ \sup_{t \le T} \Skalprod{\tilde Z_B^{N_k}}{f}^4]
\le c(\gamma,f) \tilde E[ (\int_0^T \Skalprod{\tilde X^{N_k}(s)}{1}
d s)^2]+c(X^{N_k}(0),f) \le c
$$
uniformly in $k \in \N.$ \cite{Ethier}, Chapter 7, Problem 7
yields the second claim.

Moreover, $M$ is continuous because
$\Vert \delta \tilde Z^{N_k}_B(t) \Vert_{H^{\alpha_2}(0,1)} \le
\frac{c}{N_k l_k} \to 0$ for $k \to \infty.$

The representation of the quadratic variation process of
$\Skalprod{M(t)}{f}$ follows from (\ref{martingalq}).
\end{Bew}

\begin{Lemma}\label{S19a}
For the quantities $M$ and $\tilde Y_B$ defined in Lemma
\ref{S16c}, $\tilde P$-a.s.,
$$
\tilde Y_B(t) = \int_0^t e^{(t-s)\Delta} d M(s)
$$
holds in $C(0,T,H^{\alpha_1}(0,1)).$ $M$ can be
represented as
$M(t)=\int_0^t \sqrt{\gamma\psi(s)} d W(s)$ in the
sense that $\Skalprod{M(t)}{\varphi} = \int_0^t \int_0^1
\sqrt{\gamma\psi(s,x)}\varphi(x) d W(s,x)$ for all $\varphi\in
C^\infty_{per}(0,1)$
where $W$ is a certain space-time-white noise on a possibly again
extended probability space. See also Chapter 2 of \cite{Walsh} for
an introduction to integration w.r.t. martingale measures.

Moreover, $\tilde u \in C(0,T,L^2(0,1))$ and
(\ref{darst1}) holds in this space.
\end{Lemma}

\begin{Bew}
From (\ref{partint}) we infer that
$$
\tilde Y_B^{N_k}(t) = \tilde Z_B^{N_k}(t) + \int_0^t
e^{(t-s)\Delta_{N_k}} \Delta_{N_k} \tilde Z_B^{N_k} (s) d s,
$$
and similar to the proof of Lemma \ref{S16d}
the right hand side converges to
$M(t) + \int_0^t e^{(t-s)\Delta} \Delta M(s) d s$ in
$D(0,T,H^{\alpha_2}(0,1)) + L^1(0,T,H^{\alpha_2-2}(0,1)).$
By Lemma \ref{S18}, the equality
$\tilde Y_B(t) = M(t) + \int_0^t e^{(t-s)\Delta} \Delta M(s) d s$
holds in
$D(0,T,H^{\alpha_2}(0,1)).$
Note that the stochastic integral
$\int_0^t e^{(t-s)\Delta} d M(s)$ is well defined and has a
version in $C(0,T,H^{\alpha_2}(0,1))$ according to
\cite{Kotelenez82b}. For all $m\in\Z$ we have
\begin{eqnarray*}
\Skalprod{\tilde Y_B(t)}{\varphi_m}
&=& \Skalprod{M(t)}{\varphi_m}
+\int_0^t e^{(t-s)\lambda_m} \lambda_m \Skalprod{M(s)}{\varphi_m}
d s
\\
&=&
\int_0^t e^{(t-s)\lambda_m} d\Skalprod{M(s)}{\varphi_m} =
\Skalprod{\int_0^t e^{(t-s)\Delta}d M(s)}{\varphi_m}
\end{eqnarray*}
in $C(0,T,\R)$ by a stochastic partial integration formula
whence \\
$\tilde Y_B(t) = \int_0^t e^{(t-s)\Delta} d M(s) \in
C(0,T,H^{\alpha_2}(0,1)).$ Similarly to the proof of Lemma
\ref{S8b} we have that $P_n \tilde Y_B$
tends to $\tilde Y_B$ in $D(0,T,H^{\alpha_1}(0,1))$ and
therefore
$\tilde Y_B \in C(0,T,H^{\alpha_1}(0,1)),$ $\tilde P$-a.s. This
yields the first part of the claim.

The representation of the martingale $M$ by a stochastic integral
follows by \cite{Konno}.

Since $\tilde u^{N_k}$ is bounded in $L^2(0,T,H^1_N),$ it is
bounded by Lemma \ref{HaN} in $L^2(0,T,H^\alpha)$ for all
$\alpha < \frac12.$ For possibly a subsubsequence this yields
$\tilde u^{N_k} \to \tilde u$ in $L^2(0,T,H^\alpha)$ and a.s. in
$[0,T]$ in $H^\alpha$ for all $\alpha < \frac12.$ From Lemma
\ref{S6} we infer $\tilde u \in L^\infty(0,T,L^2)$ and the claim
follows by the method of \cite{Prato}.
\end{Bew}

\section{Auxiliary Results}\label{drei}

\begin{Lemma}\label{proj}
Let $a_{0,N} = 1$ and $b_{0,N} = 0$ and
\begin{equation}\label{abnN}
a_{n,N} = \frac{N}{2\pi n} \sin(\frac{2\pi n}{N}), \,\,
b_{n,N} = \frac{N}{2\pi n} (\cos(\frac{2\pi n}{N})-1),
\end{equation}
for all $n \in \Z.$ Then we obtain for
$n=-\frac{N-1}{2},\ldots,\frac{N-1}{2}$ and $m \in \Z$
\begin{eqnarray*}
&&
\Skalprod{\varphi_{n,N}}{\varphi_m} =a_{m,N} \mbox{ for }
m=\pm n + z N,z\in\Z,m\le 0,n\le 0,
\\
&&
\Skalprod{\varphi_{n,N}}{\varphi_m} =\pm a_{m,N} \mbox{ for }
m=\pm n + z N,z\in\Z,m>0,n>0,
\\
&&
\Skalprod{\varphi_{n,N}}{\varphi_m} =-b_{m,N} \mbox{ for }
m=\pm n + z N,z\in\Z,m>0,n\le 0,
\\
&&
\Skalprod{\varphi_{n,N}}{\varphi_m} =\pm b_{m,N} \mbox{ for }
m=\pm n + z N,z\in\Z,m<0,n>0
\end{eqnarray*}
Otherwise, this scalar product is zero.
\end{Lemma}

\begin{Bew}
Elementary calculations.
\end{Bew}

\begin{Lemma}\label{HaN}
For $\alpha<\frac12$ and $f \in H_N^\alpha(0,1),$
$\Vert f \Vert^2_{H^\alpha(0,1)}
\le c \Vert f \Vert^2_{H_N^\alpha(0,1)}$
holds with a constant $c$ independent of $N$ and $f.$
\end{Lemma}

\begin{Bew}
Since this is not proved in the references we know of, we
sketch the proof. By Definition \ref{spaces}
and Lemma \ref{proj}
$$
\Vert f \Vert^2_{H^\alpha}
\le
\Skalprod{f}{1}^2+c
\sum_{n=1}^{\frac{N-1}{2}}
\Skalprod{f}{\varphi_{n,N}}^2
\Big( \sum_{l \in \Z} (a^2_{\pm n + lN,N}+b^2_{\pm n + lN,N})
(1-\lambda_{\pm n + lN})^\alpha \Big)
$$
where
$$
\sum_{l \in \Z} (a^2_{\pm n + lN,N}+b^2_{\pm n + lN,N})
(1-\lambda_{\pm n + lN})^\alpha
\le c N^2 (1-\cos(\frac{2\pi n}{N})) \sum_{l \in \Z}(\pm
n+lN)^{2(\alpha-1)}.
$$
This can easily
be estimated by $\frac{c N^2}{(1-2\alpha) n^2}
(1-\cos(\frac{2\pi n}{N})) n^{2\alpha} \le c
\frac{n^{2\alpha}}{1-2\alpha}$ which implies the claim.
\end{Bew}

\begin{Lemma}\label{S7}
For $f \in H^1_N(0,1)$ we obtain with constants not depending
on $N$ and $f$
$$
\Vert f \Vert_{H^{\frac14}_N(0,1)}
\le c \Vert f \Vert_{L^2(0,1)}^{\frac34}
\Vert f \Vert_{H^1_N(0,1)}^{\frac14}
$$
and
$$
\Vert f \Vert_{H^1_N(0,1)}
\le c \Vert f \Vert_{L^2(0,1)}
+ c \Vert \nabla^-_N f \Vert_{L^2(0,1)}
$$
\end{Lemma}

\begin{Bew}
The first inequality is an application of H\"older's inequality.
For the second, we compute for $k=1,\ldots,N$ and $m\not=0$
$$
\nabla_N^+ \varphi_{m,N}(\frac{k}{N})
= 2\pi m (a_{m,N}\varphi_{-m,N}+b_{m,N}\varphi_{m,N})
(\frac{k}{N})
$$
This yields for $m>0$
\begin{eqnarray*}
&&
\Skalprod{f}{\nabla^+_N\varphi_{m,N}}^2
+\Skalprod{f}{\nabla^+_N\varphi_{-m,N}}^2
\\
&&
=
(2\pi m)^2 \Big(\Skalprod{f}{\varphi_{m,N}}^2
+\Skalprod{f}{\varphi_{-m,N}}^2\Big)
\Big(a_{m,N}^2+b_{m,N}^2\Big).
\end{eqnarray*}
By $a_{m,N}^2+b_{m,N}^2=(-\beta_{m,N})/(2\pi m)^2$ we easily
deduce the assertion from Definition \ref{spaces}.
\end{Bew}

We now show a compactness criterion which is in some sense the
discrete equivalent to Theorem IV.4.1 in \cite{Vishik}.

\begin{Lemma}\label{S11}
Let $(u^N)_{N\in\N}$ be relatively compact in
$C(0,T,H^{\beta}(0,1))$ for some $\beta \le 0$ and
$(\nabla_N^+ u^N)_{N\in\N}$ be bounded in $L^2(0,T,L^2(0,1)).$
Then $(u^N)_{N\in\N}$ is relatively compact in $L^2(0,T,L^2(0,1)).$
\end{Lemma}

\begin{Bew}
Let $\hat u^N$ be the piecewise linear function that coincides with
$u^N$ at the points $k/N,$ $k=0,\ldots,N-1.$ Then
$\Vert\hat u^N-u^N\Vert^2_{L^2}
=\frac{1}{3 N^2} \Vert\nabla^+_N u^N\Vert^2_{L^2}$ and therefore
$\Vert \hat u^N\Vert^2_{H^1} \le c \Vert u^N \Vert^2_{L^2} +
c (1+\frac{1}{N^2})\Vert\nabla^+_N u^N\Vert^2_{L^2}$ and
$\Vert\hat u^{N+M}-\hat u^N\Vert^2_{H^\beta} \le
c \Vert u^{N+M}- u^N\Vert^2_{H^\beta}
+ \frac{c}{N^2} \Vert\nabla^+_{N+M} u^{N+M}\Vert^2_{L^2}
+ \frac{c}{N^2} \Vert\nabla^+_N u^N\Vert^2_{L^2}$.
Hence we deduce from a classical interpolation inequality that
$\forall \epsilon>0 : \exists C_\epsilon>0 : \exists N_0 \in\N:
\forall N\in\N:N>N_0 : \forall M\in\N:M>0 : $
\begin{eqnarray*}
&&
\Vert u^{N+M}- u^N\Vert^2_{L^2(0,T,L^2)}
\le \epsilon
\Big(\Vert \nabla_{N+M}^+ u^{N+M}\Vert^2_{L^2(0,T,L^2)}
+\Vert\nabla^+_N u^N\Vert^2_{L^2(0,T,L^2)} \Big)
\\
&&
\quad\quad\quad\quad\quad\quad\quad\quad\quad\quad\quad\quad
+ C_\epsilon
\Vert u^{N+M}- u^N\Vert^2_{L^2(0,T,H^\beta)}.
\end{eqnarray*}
\end{Bew}

\begin{Lemma}\label{S12}
Let $(R_N h^N)(t) = \int_0^t e^{(t-s)\Delta_N} h^N(s) d s$ for $h^N\in
H^N.$ Then for $p>1, \gamma>0$ such that $1>\frac{1}{p} + \gamma$
and for $\beta\in\R$ such that $\beta+2\gamma<0,$ $(R_N
h^N)_{N\in\N}$ is relatively compact in
$C(0,T,H^{\beta+2\gamma}(0,1))$ if $(h^N)_{N\in\N}$ is bounded in
$L^p(0,T,H^\beta(0,1)).$
\end{Lemma}

\begin{Bew}
Here we adapt the method of \cite{Gatarek}.
We obtain by a standard argument that
$\Vert e^{t \Delta_N} f\Vert_{H^{2\gamma+\beta}} \le
(1+\frac{c}{t^\gamma})\Vert f\Vert_{H^\beta}$
for $\beta\in\R$ and $\gamma>0$ with a constant $c$ not dependent
on $N$ and $f.$ An application of H\"olders
inequality then yields for sufficiently small $\tilde\epsilon >0$
that for all $t\le T$
$$
\Vert R_N h^N(t)\Vert_{H^{2\gamma+\beta+\tilde\epsilon}}
\le c(T,\tilde\epsilon,\gamma,p) \Vert h^N \Vert_{L^p(0,T,H^\beta)}
$$
whence $(R_N h^N(t))_{N\in\N}$ is relatively compact in
$H^{2\gamma+\beta}(0,1)$ for all $t\le T.$

Similarly we estimate for $s<t$
$$
\Vert \int_s^t e^{(t-\tau) \Delta_N}h^N(\tau)
d \tau\Vert_{H^{2\gamma+\beta}} \le c(T,\tilde\epsilon,\gamma,p)
\Vert h^N \Vert_{L^p(0,T,H^\beta)} |t-s|^{1-\frac{1}{p}-\gamma}.
$$
Because due to $\beta+2\gamma<0,$ the norms of
$H^{\beta+2\gamma}(0,1)$ and $H^{\beta+2\gamma}_N(0,1)$ are
equivalent, see \cite{Blount96}, we obtain
$\Vert \Delta_N e^{t \Delta_N} f\Vert_{H^{2\gamma+\beta}} \le
(1+\frac{c}{t^{\gamma+1}})\Vert f\Vert_{H^\beta}$ and
\begin{eqnarray*}
&&
\Vert\int_0^s
\Big(e^{(t-\tau) \Delta_N}-e^{(s-\tau) \Delta_N}\Big)
h^N(\tau) d \tau\Vert_{H^{2\gamma+\beta}}
\\
&&
=\Vert\int_0^s\int_{s-\tau}^{t-\tau} \Delta_N e^{\Delta_N \rho}
d\rho h^N(\tau) d \tau\Vert_{H^{2\gamma+\beta}}
\\
&&
\le c \int_0^s \Vert h^N(\tau)\Vert_{H^{\beta}}
\Big( (t-s) + (s-\tau)^{-\gamma} - (t-\tau)^{-\gamma} \Big) d \tau
\le
\end{eqnarray*}
by a technique used in \cite{Prato}, Appendix A,
$$
\le c(T,\tilde\epsilon,\gamma,p)
\Vert h^N \Vert_{L^p(0,T,H^\beta)} |t-s|^{1-\frac{1}{p}-\gamma}
$$
for $|t-s| \le 1.$ This shows that $(R_N h^N)_{N\in\N}$ is
equicontinuous in $H^{2\gamma+\beta}(0,1)$ on $[0,T].$
\end{Bew}

\begin{Lemma}\label{d2p}
Convergence in $D(0,T,H^\alpha(0,1))$ implies convergence in \\
$L^p(0,T,H^\alpha(0,1)),$ $p \ge 1.$
\end{Lemma}

\begin{Bew}
First, $D(0,T,H^\alpha)\subset L^\infty(0,T,H^\alpha)$
algebraically due to the existence of left limits and right
continuity in $D(0,T,H^\alpha).$ Let $(f_n)_{n\in\N}$ be a
sequence in $D(0,T,H^\alpha)$ which converges to $f$ in that
space. Then according to \cite{Ethier} there exists a sequence of
strictly increasing Lipschitz continuous functions
$(\rho_n)_{n\in\N}$ with $\rho_n(0)=0$ and $\rho_n(T)=T$ such that
$\lim_{n\to\infty} \sup_{t\le T} |\rho_n(t)-t|=0$ and
$\lim_{n\to\infty} \sup_{t\le T} \Vert f_n(t)-f(\rho_n(t))
\Vert_{H^\alpha} = 0.$ Because of 5.5.1 in \cite{Ethier},
$f(\rho_n(t))\to f(t)$ a.e. for $n\to\infty$ and therefore
$\int_0^T\Vert f(t)-f(\rho_n(t)) \Vert^p_{H^\alpha} d t$ due to
the integrable bound $c\Vert f\Vert^p_{L^\infty(0,T,H^\alpha)}.$
\end{Bew}


\begin{thebibliography}{50}

\bibitem{Adams}
R. Adams.
{\em Sobolev Spaces.} Volume 65 of {\em Pure and Applied Mathematics}.
Academic Press Inc.,1978.

\bibitem{Arnold}
L. Arnold, M. Theodosopulu.
Deterministic limit of the stochastic model of chemical
reactions with diffusion.
{\em Advances in Applied Probability} 12 (1980), 367--379.

\bibitem{Bertini}
L. Bertini, G. Giacomin.
Stochastic Burgers and KPZ equations from particle systems.
{\em Communications in Mathematical Physics} 183 (1997), 571--607.

\bibitem{Biler}
P. Biler, T. Funaki, and W. A. Woyczynski.
Fractal Burgers equations.
{\em Journal of Differential Equations} 148 (1998), 9--46.

\bibitem{Blountdiss}
D. Blount.
{\em Comparison of a Stochastic Model of a Chemical
Reaction with Diffusion and the Deterministic Model.}
PhD thesis, University of Wisconsin-Madison, 1987.

\bibitem{Blount91}
D. Blount.
Comparison of stochastic and deterministic models of a
linear chemical reaction with diffusion.
{\em Annals of Probability} 19 (1991), no. 4, 1440--1462.

\bibitem{Blount92}
D. Blount.
Law of large numbers in the supremum norm for a chemical
reaction with diffusion.
{\em Annals of Applied Probability} 2 (1992), no. 1,
131--141.

\bibitem{Blount93}
D. Blount.
Limit theorems for a sequence of nonlinear
reaction-diffusion systems.
{\em Stochastic Processes and their Applications}
45 (1993), 193--207.

\bibitem{Blount94}
D. Blount.
Density-dependent limits for a nonlinear reaction-diffusion
model.
{\em The Annals of Probability} 22 (1994), no. 4,
2040--2070.

\bibitem{Blount96}
D. Blount.
Diffusion limits for a nonlinear density dependent
space-time population model.
{\em The Annals of Probability} 24 (1996), no. 2,
639--659.

\bibitem{Bonnet}
G. Bonnet, R. J. Adler.
The Burgers superprocess.
{\em Preprint} (2002).

\bibitem{Calderoni}
P. Calderoni, M. Pulvirenti.
Propagation of chaos for Burgers equation.
{\em Ann. Inst. Henri Poincare, Nouv. Ser.
Sect A} 39 (1983), 85--97.

\bibitem{Cap}
M. Capinski, D. Gatarek.
Stochastic equations in Hilbert space with applications
to Navier-Stokes equation in any dimensions.
{\em Journal of Functional Analysis}, 126 (1994),
no. 1, 26--35.


\bibitem{Prato}
G. Da Prato, J. Zabczyk.
{\em Ergodicity in Infinite Dimensional Systems.}
Cambridge Univ. Press, 1996.

\bibitem{Daw}
D. Dawson.
{\em Measure Valued Markov Processes.}
Lecture Notes in Mathematics 1541, Springer-Verlag,
1993, 1--260.

\bibitem{DeMasi89}
A. De Masi, P. A. Ferrari, M. E. Vares.
A microscopic model of interface related to the Burgers
equation.
{\em J. Stat. Phys.} 55 (1989), no. 3/4, 601--609

\bibitem{DeMasi}
A. De Masi, N. Ianiro, A. Pellegrinotti, and E. Presutti.
{\em A Survey of the Hydrodynamical Behaviour
of Many Particle Systems.}
Studies in Statistical Mechanics, North Holland, 1984,
123--294.

\bibitem{Ethier}
S. N. Ethier, T. G. Kurtz.
{\em Markov Processes. Characterization and Convergence.}
John Wiley \& Sons, 1986.

\bibitem{Gatarek}
D. Gatarek.
A note on nonlinear stochastic equations in Hilbert spaces.
{\em Statistics \& Probability Letters} 17 (1993), 387--394.

\bibitem{Grecksch}
W. Grecksch, B. Schmalfu\ss.
Approximation of the stochastic Navier-Stokes equation.
{\em Comp. Appl. Math.} 15 (1996), no. 3, 227--239.

\bibitem{Gugg}
C. Gugg.
{\em Approximation of Stochastic Partial Differential
Equations and Turbulence in Fluids.}
PhD Thesis,
Wi\ss ner-Verlag, Augsburg, 2001.

\bibitem{Kielhoefer}
C. Gugg, H. Kielh\"ofer, and M. Niggemann.
On the approximation of the stochastic Burgers equation.
To appear in {\em Commun. Math. Phys.}


\bibitem{Kipnis}
C. Kipnis, S. Olla, and S. Varadhan.
Hydrodynamics and large deviation for simple exclusion
processes.
{\em Communications on Pure and Applied Mathematics} 42 (1989),
115--137.

\bibitem{Konno}
N. Konno, T. Shiga.
Stochastic partial differential equations for some
measure-valued diffusions.
{\em Probability Theory and Related Fields} 79 (1988),
201--225.

\bibitem{Kotelenez82}
P. Kotelenez.
{\em Law of large numbers and central limit theorem
for linear chemical reactions with diffusion.}
PhD Thesis,
Universit\"at Bremen, 1982.

\bibitem{Kotelenez82b}
P. Kotelenez.
A submartingale type inequality with applications to stochastic
evolution equations.
{\em Stochastics} 8 (1982), 139--151.

\bibitem{Kotelenez86a}
P. Kotelenez.
Law of large numbers and central limit theorem for linear chemical
reactions with diffusion.
{\em The Annals of Probability} 14 (1986), no. 1, 173--193.

\bibitem{Kotelenez86b}
P. Kotelenez.
Linear parabolic differential equations as limits of
space-time jump Markov processes.
{\em J. Math. Anal. Appl.} 116 (1986), 42--76.

\bibitem{Kotelenez88}
P. Kotelenez.
High density limit theorems for nonlinear chemical reactions with
diffusion.
{\em Probability Theory and Related Fields} 78 (1988), 11--37.

\bibitem{Kouritzin02}
M. Kouritzin, H. Long.
Convergence of Markov chain approximations to stochastic reaction
diffusion equations.
{\em Preprint} (2002).


\bibitem{Liggett}
T. Liggett.
Stochastic models of interacting systems.
{\em Ann. Probability} 25 (1997), 1--29.

\bibitem{Meleard}
S. M\'el\'eard.
Monte-Carlo approximations for 2d Navier-Stokes equations
with measure initial data.
{\em Probability Theory and Related Fields} 121 (2001),
no. 3, 367--388.

\bibitem{Oel85}
K. Oelschl\"ager.
A law of large numbers for moderately interacting diffusion
processes.
{\em Z. Wahrscheinlichkeitstheorie verw. Gebiete} 69 (1985), 279--322.

\bibitem{Vishik}
M. Vishik, A. Fursikov.
{\em Mathematical Problems of Statistical Hydromechanics.}
Kluwer, 1988.

\bibitem{Walsh}
J. B. Walsh.
{\em An Introduction to Stochastic Partial Differential
Equations.}
Lecture Notes in Mathematics 1180, Springer-Verlag, 1986,
265--437.

\end{thebibliography}
\end{document}